\documentclass[a4paper]{amsart}

\usepackage{times}

\usepackage[T1]{fontenc}
\usepackage[latin1]{inputenc}
\usepackage{hyperref}
\usepackage[arrow,matrix,curve]{xy}

\def\factor#1.#2.{\left. \raise 2pt\hbox{$#1$} \right/\hskip -2pt\raise -2pt\hbox{$#2$}}

\theoremstyle{plain}
\newtheorem{thm}{Theorem}[section]
\numberwithin{equation}{section} 
\numberwithin{figure}{section} 
\theoremstyle{plain}
\newtheorem{cor}[thm]{Corollary} 
\newtheorem{lem}[thm]{Lemma} 
\theoremstyle{plain}
\newtheorem{prop}[thm]{Proposition} 
\newtheorem{fact}[thm]{Fact} 
\theoremstyle{remark}
\newtheorem{rem}[thm]{Remark}
\theoremstyle{remark}
\newtheorem{notation}[thm]{Notation}

\newcommand{\NZ}{\mbox{$\mathbb{N}$}}

\DeclareMathOperator{\Aut}{Aut}
\DeclareMathOperator{\codim}{codim}
\DeclareMathOperator{\Hom}{Hom}
\DeclareMathOperator{\Image}{Image}
\DeclareMathOperator{\Pic}{Pic}
\DeclareMathOperator{\Spec}{Spec}

\newcommand\sE{{\mathcal E}}
\newcommand\sF{{\mathcal F}}
\newcommand\sO{{\mathcal O}}

\title {Holomorphic maps onto varieties of non-negative Kodaira
  dimension}

\author{Jun-Muk Hwang, Stefan Kebekus, Thomas Peternell}

\date{\today}

\address{Jun-Muk Hwang, Korea Institute for Advanced Study, 207-43
  Cheongnyangni-dong, Seoul, 130-722, Korea}
\email{jmhwang@kias.re.kr}

\address{Stefan Kebekus, Mathematisches Institut, Universit\"at zu
K\"oln, Weyertal 86-90, 50931 K\"oln, Germany}
\email{stefan.kebekus@math.uni-koeln.de}
\urladdr{http://www.mi.uni-koeln.de/\~{}kebekus}

\address{Thomas Peternell, Institut für Mathematik, Universität
  Bayreuth, 95440~Bayreuth, Germany}
\email{thomas.peternell@uni-bayreuth.de}

\thanks{Jun-Muk Hwang was supported by the Korea Research Foundation
  Grant (KRF-2002-070-C00003). Stefan Kebekus was supported by a
  Heisenberg-Fellowship of the DFG. Jun-Muk Hwang, Stefan Kebekus and
  Thomas Peternell were supported in part by the program ``Globale
  Methoden in der Komplexen Analysis'' of the DFG}

\begin{document}

\maketitle
\tableofcontents

\section{Introduction and Statement of Results}

A classical result in complex geometry says that the automorphism
group of a manifold of general type is discrete \cite{Mat63}. It is
more generally true that there are only finitely many surjective
morphisms between two fixed projective manifolds $X$ and $Y$ of
general type \cite{KO75}.

Rigidity of surjective morphisms, and the failure of a morphism to be
rigid have been studied by a number of authors, the most general
results being those of Borel and Narasimhan \cite{BN67}. For target
manifolds $Y$ with Chern numbers $c_1(Y) = 0$ and $c_{\dim Y}(Y) \ne
0$, rigidity has been shown by Kalka, Shiffman and Wong \cite{KSW81}.
These results have recently been generalized by Hwang \cite{Hwa03} to
the case where $Y$ is a compact Kähler manifold with $c_1(Y) = 0$.
Although in Hwang's setup deformations need not be rigid, he is able
to give a good description of the space of surjective morphisms.

In this paper we give a complete description of the space of
surjective morphisms in the general setup where $Y$ is a normal
projective variety that is not covered by rational curves. Our main
result, Theorem~\ref{thm:main}, states that surjective morphisms are
rigid, unless there is a clear geometric reason for it.
\begin{itemize}
\item Deformations of surjective morphisms between normal projective
  varieties are unobstructed unless the target variety is covered by
  rational curves.

\item If the target is not covered by rational curves, then surjective
  morphisms are infinitesimally rigid, except for those morphisms that
  factor via a variety with positive-dimensional automorphism group.
\end{itemize}

\begin{notation}
  If $X$ and $Y$ are normal compact complex varieties, $\Hom(X,Y)$
  denotes the space of holomorphic maps $X \to Y$ and $\Hom_s(X,Y)$
  the space of surjective holomorphic maps. Given a morphism $f \in
  \Hom(X,Y)$, let $\Hom_f(X,Y)$ be the connected component of
  $\Hom(X,Y)$ that contains $f$.
\end{notation}

\begin{thm}\label{thm:main}
  Let $X$ be a normal compact complex variety and $Y$ be a projective
  normal variety which is not covered by rational curves. If $f: X \to
  Y$ is a surjective morphism, then there exists a factorization
  $$
  \xymatrix{
    X  \ar[r]_{\alpha} \ar@/^0.3cm/[rr]^{f} &
    Z  \ar[r]_{\beta} &
    Y
  }
  $$
  where
  \begin{enumerate}
  \item $\beta$ is a finite morphism which is étale outside of the
    singular set of $Y$

  \item if $\Aut^0(Z)$ is the maximal connected subgroup of the
    automorphism group of $Z$, then $\Aut^0(Z)$ is an Abelian variety,
    and the natural morphism
    $$
    \factor \Aut^0(Z) . \Aut(Z/Y) \cap \Aut^0(Z) . \to \Hom_f(X,Y)
    $$
    is isomorphic.
  \end{enumerate}
  In particular, all deformations of surjective morphisms $X \to Y$ are
  unobstructed, and the associated components of $\Hom(X,Y)$ are smooth
  Abelian varieties.
\end{thm}

\begin{cor}\label{cor:main}
  In the setup of Theorem~\ref{thm:main}, if $Y$ is smooth, then $Y$
  has a finite étale cover of the form $T \times W$, where $T$ is an
  Abelian variety of dimension $h^0\bigl(X,f^*(T_Y)\bigr)$ and
  $$
  \dim \Hom_f(X,Y) \leq \dim Y - \kappa(Y),
  $$
  where $\kappa(Y)$ is the Kodaira dimension.
\end{cor}

\begin{rem}
  We conjecture that Theorem~\ref{thm:main} and
  Corollary~\ref{cor:main} are true when $Y$ is a compact Kähler
  manifold of nonnegative Kodaira dimension. Our proof needs the
  projectivity assumption because it employs Miyaoka's
  charactierzation of unirueledness.
\end{rem}

The following corollaries are immediate consequences of
Corollary~\ref{cor:main}.

\begin{cor}
  Let $Y$ be a projective manifold which is not uniruled. If $\pi_1(Y)
  $ is finite, then for each connected normal compact complex variety
  $X$ the space $\Hom_s(X,Y)$ is discrete .
\end{cor}

\begin{cor}
  Let $Y$ be a projective $n-$dimensional manifold which is not
  uniruled. If $c_n(Y) \ne 0,$ then for each connected normal compact
  complex variety $X$ the space $\Hom_s(X,Y)$ is discrete.
\end{cor}

\section{Proof of Theorem~\ref{thm:main}}

\subsection{Step 1: Setup}

Let $X$ be a normal variety. Then the tangent sheaf $T_X$ is by
definition the dual of the sheaf $\Omega_X^1$ of differentials.
If $f: X \to Y$ is holomorphic, we consider $\Hom_f(X,Y)$, the
connected component of $\Hom(X,Y)$ that contains $f$. If $f$ is
additionally surjective, since $X$ is reduced, it is then
well-known that
$$
T_{\Hom_f(X,Y)}|_f \simeq \Hom\bigl(f^*(\Omega^1_Y), \sO_X\bigr).
$$
See e.g.~\cite[I, Thm.~2.16]{K96} for a proof in the algebraic case.
We note that if $Y$ is smooth, then $\Hom\bigl( f^*(\Omega^1_Y),
\sO_X\bigr) \cong H^0\bigl(X, f^*(T_Y)\bigr)$.

If in the set-up of Theorem~\ref{thm:main}, there are no infinitesimal
deformations of the morphism $f$, i.e.~if $\Hom\bigl( f^*(\Omega^1_Y),
\sO_X\bigr) = \{0\}$, there is nothing to prove. We will therefore
assume throughout that $\Hom\bigl( f^*(\Omega^1_Y), \sO_X\bigr) \not =
\{0\}$.

\subsection{Step 2: Reduction to a finite morphism}

In this section we reduce the proof of Theorem~\ref{thm:main} to the
case that the morphism $f$ is finite. To this end, we will consider
the Stein factorization of $f$,
\begin{equation}
  \label{eq:SteinFact}
\xymatrix{
  X  \ar[rr]_{g \text{, conn. fibers}} \ar@/^0.3cm/[rrrr]^{f} &&
  W  \ar[rr]_{h \text{, finite}} &&
  Y,
}
\end{equation}
assume that Theorem~\ref{thm:main} holds for the finite morphism $h$,
and show that $\Hom_f(X,Y)$ and $\Hom_h(W,Y)$ are then naturally
isomorphic. The argumentation is based on the following elementary
observation whose proof we leave to the reader.

\begin{fact}\label{fact:criterionForEtale}
  Let $h: S \to B$ be a morphism of complex spaces. Assume that $S$ is
  smooth and compact, $B$ is connected and that the associated
  morphism between the Zariski tangent spaces is everywhere
  isomorphic. Then $h$ is surjective and étale.

  In particular, $h$ is an isomorphism if it is injective.
\end{fact}

In order to apply Fact~\ref{fact:criterionForEtale}, observe
that the Stein factorization~\eqref{eq:SteinFact} yields a canonical
morphism of complex spaces
$$
\begin{array}{rccc}
A: & \Hom_h(W,Y) & \to     & \Hom_f(X,Y) \\
   & \gamma      & \mapsto & \gamma \circ g
\end{array}
$$
which is injective because $g$ is surjective. If $\gamma \in
\Hom_h(W,Y)$ is any morphism, it is known that associated morphism
between the Zariski tangent spaces at $\gamma$ and $\gamma \circ g$
$$
TA:
\underbrace{T_{\Hom(W,Y)}|_{\gamma}}_{\Hom\bigl(\gamma^*(\Omega^1_Y),
  \sO_W\bigr)} \to \underbrace{T_{\Hom(X,Y)}|_{\gamma \circ
    g}}_{\Hom\bigl(g^*\gamma^*(\Omega^1_Y),\sO_X\bigr)}
$$
is the pull-back via $g$. Since $g$ has connected fibers, $g_*(\sO_X)
= \sO_W$, and since $g_*$ and $g^*$ are adjoint functors,
\cite[p.~110]{Ha77}, this map is isomorphic.

If Theorem~\ref{thm:main} holds for the finite morphism $h$,
$\Hom_h(W,Y)$ will be a projective manifold. By
Fact~\ref{fact:criterionForEtale}, the morphism $A$ will then be
isomorphic, and Theorem~\ref{thm:main} will hold for $f$, too.  We
are therefore reduced to showing Theorem~\ref{thm:main} under the
additional assumption that $f$ is finite. We maintain this assumption
throughout the rest of the proof.

\begin{rem}
  If $f$ is finite and $H \in \Pic(Y)$ ample, then $f^*(H)$ will again
  be ample. Thus, the assumption that $f$ is finite implies that $X$
  is projective. We can therefore argue in the algebraic category for
  the remainder of the proof.
\end{rem}

\subsection{Step 3: Further setup}

In the sequel we will use the following notation:

\begin{align*}
  \text{ample bundle:} &  &
  H   & \ldots \text{an ample line bundle on $Y$} \\
  \text{exceptional sets:} &  &
      X_s & \ldots \text{singular locus of $X$} \\
  & & Y_s & := f(X_s) \cup \{ \text{singular locus of $Y$} \}\\
  \text{open sets:} &  &
  Y_0 & := Y \setminus Y_s \\
  & & X_0 & := f^{-1}(Y_0) \\
  & & f_0 & := f|_{X_0} : X_0 \to Y_0
\end{align*}

It is well-known that the finite morphism $f_0$ defines a vector
bundle on the quasi-projective target manifold $Y_0$.

\begin{fact}
  The trace map gives a splitting
  $$
  (f_0)_* \sO_{X_0} \cong \sO_{Y_0} \oplus \sE_0^*
  $$
  where $\sE_0^*$ is a vector bundle on $Y_0$. In particular, the
  projection formula gives
  $$
  (f_0)_* (f_0)^* T_{Y_0} \cong T_{Y_0} \oplus (\sE_0^* \otimes T_{Y_0}).
  $$
\end{fact}

\begin{rem}
  The exceptional set $Y_s$ is of codimension $\geq 2$. Thus, if $m
  \in \NZ$ is sufficiently large and $H_1, \ldots, H_{\dim Y-1} \in
  |mH|$ are general members, then the general complete intersection
  curve
  $$
  C := H_1 \cap \ldots \cap H_{\dim Y -1}
  $$
  does not intersect $Y_s$. In particular, the vector bundle
  $\sE_0^*$ is defined all along $C$.
\end{rem}

\subsection{Step 4: Construction of the étale cover}

In this section we construct a factorization of the morphism $f$,
which we assume to be finite, via an étale cover of $Y$. The important
properties of the construction are summarized in the following
proposition.

\begin{prop}\label{prop:constrCover}
  In the setup of Theorem~\ref{thm:main}, there exists a canonical
  factorization of $f$ via a finite morphism $\beta$ that is étale
  outside of the singular set of $Y$,
  $$
  \xymatrix{
    X  \ar[r]_{\alpha} \ar@/^0.3cm/[rr]^{f} &
    Z  \ar[r]_{\beta} &
    Y
  }
  $$
  such that all infinitesimal deformations of $f$ come from pull-backs
  of vector fields on $Z$, i.e.~that the natural injective morphism
  $$
  \underbrace{\Hom\bigl(\Omega^1_Z, \sO_Z\bigr)}_{\text{vector fields}} \to \Hom \bigl( \Omega^1_Z,
  \alpha_* (\sO_X) \bigr) \cong \Hom \bigl(\alpha^* (\Omega^1_Z), \sO_X \bigr) \cong
  \underbrace{\Hom \bigl( f^*(\Omega^1_Y), \sO_X \bigr)}_{\text{infinit.~deformations}}
    $$
    is isomorphic.
\end{prop}

\begin{rem}
  In the formulation of Proposition~\ref{prop:constrCover} we have
  identified $\Hom \bigl(\alpha^* (\Omega^1_Z), \sO_X \bigr)$ and
  $\Hom \bigl( f^*(\Omega^1_Y), \sO_X \bigr)$. For this, we use the
  assumptions that $f$ is finite and that $\beta$ is étale outside of
  a set of codimension 2: the (reflexive) sheafs $(\alpha^*
  (\Omega^1_Z))^\vee$ and $(f^*(\Omega^1_Y))^\vee$ agree in
  codimension 1. Since $X$ is normal, they must be isomorphic.

  If $Y$ is smooth, then $Z$ must also be smooth and the natural
  morphism discussed in Proposition~\ref{prop:constrCover} is simply
  the pull-back map
  $$
  \alpha^* : H^0\bigl(Z, T_Z\bigr) \to H^0\bigl(X, \alpha^*(T_Z)\bigr) \cong H^0\bigl(X, f^*(T_Y)\bigr).
  $$
\end{rem}

We start the proof of Proposition~\ref{prop:constrCover} with the
following lemma which links the existence of elements in $\Hom \bigl(
f^*(\Omega^1_Y), \sO_X \bigr)$ that do not come from vector fields to
the structure of the bundle $\sE_0$.

\begin{lem}\label{lem:Enefnotample}
  Assume that there exists an infinitesimal deformation $\sigma \in
  \Hom \bigl( f^*(\Omega^1_Y), \sO_X \bigr)$ which does \emph{not}
  come from the pull-back of a vector field on $Y$. Then, if $C$ is a
  general complete intersection curve and $\sE_0$ the dual of
  $\sE_0^*$, the restriction $\sE_0|_C$ is nef, but not ample.
\end{lem}

\begin{proof}
  Since $C$ is not contained in the branch locus, the fact that
  $\sE_0|_C$ is nef is shown in \cite[Thm. A of the appendix by
  R.~Lazarsfeld]{PS00} ---as we need only the nefness on a general
  curve, we could also use the general semi-positivity theorem of
  Viehweg for images of relative dualizing sheaves.

  Recall that $\codim_X X \setminus X_0 \geq 2$. Sections in a
  reflexive sheaf which are defined on $X_0$ therefore extend uniquely
  to all of $X$. This yields identifications
  \begin{align*}
     \Hom \bigl( f^*(\Omega^1_Y), \sO_X \bigr) & = H^0\bigl(X_0, (f_0)^*(T_{Y_0})\bigr) \\
      & = H^0\bigl(Y_0, T_{Y_0} \oplus (\sE_0^* \otimes T_{Y_0})\bigr) \\
      & = H^0\bigl(Y, T_Y\bigr) \oplus H^0\bigl(Y_0, \sE_0^* \otimes T_{Y_0}\bigr).
  \end{align*}
  Since we assume that the infinitesimal deformation $\sigma$ does not
  come from the pull-back of a vector field, we obtain a section
  $\tilde \sigma \in H^0\bigl(Y_0, \sE_0^* \otimes T_{Y_0}\bigr)$,
  i.e., a morphism of vector bundles
  $$
  \tilde \sigma : \sE_0 \to T_{Y_0}.
  $$
  After removing further sets of codimension 2, if necessary, we may
  assume without loss of generality that
  $$
  \mathcal F := \Image(\tilde \sigma) \subset
  T_{Y_0}
  $$
  is a locally free subsheaf of $T_{Y_0}$. The restriction of its dual
  to a general complete intersection curve, $\mathcal F^*|_C$, is then
  a torsion-free quotient of $\Omega^1_{Y_0}|_C$, which, by Miyaoka's
  celebrated theorem~\cite[cor.~6.4]{Miyaoka87} (see also
  Theorem~9.0.1 of Shepherd-Barron's article in
  \cite{SecondAsterisque}) has non-negative degree. Equivalently, we
  can say that $\mathcal F|_C$ has non-positive degree. But $\mathcal
  F|_C$ is a quotient of $\sE_0|_C$ and should therefore have positive
  degree if $\sE_0|_C$ was ample. We conclude that $\sE_0|_C$ is not
  ample.
\end{proof}

The existence of a factorization of $f$ via a cover of $Y_0$ now
follows from the argumentation of \cite[proof of Prop.~3.8]{PS03}. For
the reader's convenience, we reproduce the proof here.

\begin{lem}\label{lem:firstcover}
  In the setup of lemma~\ref{lem:Enefnotample}, after perhaps removing
  further subsets of codimension two, if necessary, the morphism $f_0$
  factors via an étale cover $Y_0^{(1)} \to Y_0$, which is not an
  isomorphism.
\end{lem}

\begin{proof}
  To factorize the morphism $f_0$, it suffices to find a coherent
  subsheaf $\sF \subset \sE_0^*$ such that $\sO_{Y_0} \oplus \sF
  \subset \sO_{Y_0} \oplus \sE_0^* \cong (f_0)_*\sO_{X_0}$ is a sheaf
  of $\sO_{Y_0}$-algebras, i.e.~closed under the multiplication map
  $$
  \mu : (\sO_{Y_0}\oplus \sE_0^*) \otimes (\sO_{Y_0}\oplus \sE_0^*)
  \to \sO_{Y_0}\oplus \sE_0^*
  $$
  We can then set $Y_0^{(1)} := \Spec \sF$. If $\sF \subset
  \sE_0^*$ is a sub-vectorbundle that has degree zero on the general
  complete intersection curve, then it follows that the natural
  morphism $Y_0^{(1)} \to Y_0$ is étale.

  As a first step towards the construction of $\sF$, we fix a complete
  intersection curve $C \subset Y_0$ and construct $Y_0^{(1)}$ only
  over $C$. Since the restriction $\sE_0|_C$ is nef, but not ample, it
  follows from \cite[Lem.~2.3]{PS00} that there exists a unique
  maximal ample subbundle $V_C \subset \sE_0|_C$ such that the
  quotient $\sE_0|_C / V_C$ has degree zero. Let $\sF_C \subset
  \sE_0^*$ be the kernel of the associated map $\sE_0^*|_C \to V_C^*$
  which is a sub-vectorbundle of degree zero. It is then clear that
  $\sO_C \oplus \sF_C \subset \sO_C \oplus \sE_0^*|_C$ is closed under
  multiplication, as the map
  $$
  \mu' : \underbrace{(\sO_C \oplus \sF_C) \otimes (\sO_C \oplus
    \sF_C)}_{\text{degree zero}} \to \underbrace{\factor \sO_C \oplus
    \sE_0^*|_C . \sO_C \oplus \sF_C .}_{\text{negative degree}}
  $$
  is necessarily zero.

  To end the proof of Lemma~\ref{lem:firstcover}, we need to extend
  the sub-vectorbundle $V_C \subset \sE_0|_C$ to all of $Y_0$, i.e.~we
  need to find a sub-vectorbundle $T \subset \sE_0$ such that for a
  general complete intersection curve $C' \subset Y_0$, the
  restriction $T|_{C'} \subset \sE_0|_{C'}$ is the unique maximal
  ample subbundle. For this, consider the Harder-Narasimhan filtration
  of $\sE_0|_C$,
  $$
  0 = \sE_0|_C^{(0)} \subset \sE_0|_C^{(1)} \subset \cdots \subset
  \sE_0|_C^{(\ell)} = \sE_0|_C.
  $$
  It is an elementary computation to see that there exists a number
  $k$ such that $V_C = \sE_0|_C^{(k)}$. In this setup, after removing
  further subsets of codimension two, if necessary, the theorem of
  Mehta-Ramanathan \cite[Thm.~9.1.1.7]{SecondAsterisque} (see also
  \cite{MR82}) guarantees that $V_C$ extends to all of $Y_0$, as
  required.
\end{proof}

\begin{rem}
  János Kollár pointed out to us that the proof of
  Lemma~\ref{lem:firstcover} really shows that if an antinef vector
  bundle on a curve has a section after pull back, then it has a
  section after an étale pull back
\end{rem}

It is a classical result that the cover $Y_0^{(1)} \to Y_0$ can be
extended to $Y$.

\begin{cor}\label{cor:constrOfCover}
  In the setup of Lemma~\ref{lem:Enefnotample}, the morphism $f$
  factors via a normal variety $Y^{(1)}$,
  $$
  \xymatrix{
    X  \ar[r]_{a} \ar@/^0.3cm/[rr]^{f} &
    Y^{(1)} \ar[r]_{b} &
    Y
  }
  $$
  where $b$ is a finite morphism of degree $>1$, étale outside of
  the singular locus of $Y$.
\end{cor}

\begin{proof}
  The factorization for $f_0 : X_0 \to Y_0$ via an unbranched cover
  $b_0: Y^{(1)}_0 \to Y_0$ is shown in Lemma~\ref{lem:firstcover}.
  Since $Y$ is normal, \cite[Thm.~3.5]{DG94} says that there exists a
  unique normal compactification $Y^{(1)} \supset Y^{(1)}_0$ with a
  finite morphism $b: Y^{(1)} \to Y$ that extends $b_0$ and is étale
  outside of the singular set.

  The proof is finished if we show that the associated rational map
  $a: X \dasharrow Y^{(1)}$ is a morphism. That, however, follows from
  that fact that $f = b\circ a$ is a morphism and that $b$ is finite.
\end{proof}

\begin{proof}[Proof of Proposition~\ref{prop:constrCover}]
  If all infinitesimal deformations of $f$ come from pull-backs of
  vector fields on $Z$, i.e.~if $\Hom \bigl( f^*(\Omega^1_Y), \sO_X
  \bigr) \cong H^0(Z,T_Z)$, there is nothing to prove: set $Z=Y$.

  If there exists an infinitesimal deformation $\sigma_1 \in \Hom
  \bigl( f^*(\Omega^1_Y), \sO_X \bigr)$ that is not a pull-back of a
  vector field on $Y$, apply Corollary~\ref{cor:constrOfCover} and
  obtain a factorization of $f$ via a cover $b: Y^{(1)} \to Y$. If
  there is a section $\sigma_2 \in \Hom \bigl( f^*(\Omega^1_Y), \sO_X
  \bigr) = \Hom \bigl( a^*(\Omega^1_{Y^{(1)}}), \sO_X \bigr)$, which
  is not the pull-back of a vector field on $Y^{(1)}$, apply
  Corollary~\ref{cor:constrOfCover} again to the morphism $a:X \to
  Y^{(1)}$. Proceed inductively, creating a sequence of covers $$
  \xymatrix{ X \ar[r] \ar@/^0.5cm/[rrrrr]^{f} & Y^{(d)} \ar[r] &
    Y^{(d-1)} \ar[r] & \ldots \ar[r] & Y^{(1)} \ar[r] & Y } $$
  The
  process terminates because the degree of $f$ is finite. Let $Z :=
  Y^{(d)}$ be the terminal variety.
\end{proof}

\subsection{Step 5: end of proof}

The factorization of $f$ given in Proposition~\ref{prop:constrCover}
yields a natural morphism
$$
\begin{array}{rccc}
\iota: & \Aut^0(Z) & \to & \Hom_f(X,Y) \\
       & g & \mapsto & \beta \circ g \circ \alpha
\end{array}
$$
Note that the proof of Theorem~\ref{thm:main} is finished if we
show that $\iota$ is étale. This will be guaranteed by
Fact~\ref{fact:criterionForEtale} as soon as the following two
prerequisites are verified.

\begin{description}
\item[$\boldsymbol \Aut^0(Z)$ is proper] By assumption, the variety
  $Y$ is not uniruled. Thus, the variety $Z$ is also not uniruled, and
  it follows from \cite{Ros56} that the automorphism group $\Aut^0(Z)$
  does not contain an algebraic subgroup which is isomorphic to
  $\mathbb C$ or to $\mathbb C^*$. The group $\Aut^0(Z)$ must
  therefore be an Abelian variety, $\Aut^0(Z) \cong \mathbb
  C^m/\Gamma$. In particular, $\Aut^0(Z)$ is proper.

\item[The tangent morphism is everywhere isomorphic] It is known that
  for any automorphism $g \in \Aut^0(Z)$, the morphism $T\iota$
  between Zariski tangent spaces,
  $$
  T \iota: \underbrace{T_{\Aut^0(Z)}|_g}_{H^0 \bigl(Z, g^*(T_Z) \bigr)} \to
  \underbrace{T_{\Hom(X,Y)}|_{\beta \circ g \circ \alpha}}_{\Hom\bigl(
    (\beta \circ g \circ \alpha)^*(\Omega^1_Y), \sO_X \bigr)}
  $$
  is given by the natural injective morphism of sheaves
  $$
\quad\quad  \underbrace{\Hom\bigl(g^*(\Omega^1_Z), \sO_Z\bigr)}_{H^0 \bigl(Z,
    g^*(T_Z) \bigr)} \to \Hom \bigl( g^*(\Omega^1_Z), \alpha_* (\sO_X)
  \bigr) \cong \underbrace{\Hom \bigl(\alpha^* g^*(\Omega^1_Z), \sO_X
    \bigr)}_{\Hom\bigl( (\beta \circ g \circ \alpha)^*(\Omega^1_Y),
    \sO_X \bigr)} $$
  Since $g$ is an automorphism, $g^*(\Omega^1_Z) \cong \Omega^1_Z$,
  and Proposition~\ref{prop:constrCover} asserts that this morphism is
  indeed isomorphic.
\end{description}

This ends the proof of Theorem~\ref{thm:main}. \qed

\section{Proof of Corollary~\ref{cor:main}}

In the setup of Corollary~\ref{cor:main}, the varieties $Y$ and $Z$
are smooth. If $\Aut^0(Z)$ is trivial, i.e.~$H^0(Z, T_Z)=0$, there is
nothing to prove. Otherwise, since $Z$ is not uniruled, the fact that
$Z$ is a Torus-Seifert fibration follows from \cite[Thm.~4.9]{Lib78}.
By \cite[Thm.~4.10]{Lib78}, we have that
$$
\underbrace{\dim H^0(Z,
  T_Z)}_{= \dim \Hom_f(X,Y)} + \kappa(Z) \leq \underbrace{\dim Z}_{=
  \dim Y}
$$
and therefore
$$
\dim \Hom_f(X,Y) \leq \dim Y - \kappa(Z)
\leq \dim Y - \kappa(Y).
$$
\qed

\end{document}